\documentclass[11pt]{article}
\usepackage{natbib,fullpage,amsmath,amssymb,graphicx,hyperref,pgf,pgfarrows,pgfnodes}
\newcommand{\m}[1]{\mbox{\bf{#1}} }
\newcommand{\mc}[1]{\ensuremath{{\mathcal #1}}}

\begin{document}
\title{Extending the rank likelihood for semiparametric copula estimation}
\author{Peter D. Hoff \thanks{Departments of Statistics, Biostatistics and the
Center for Statistics and the Social Sciences,
University of Washington,
Seattle, Washington 98195-4322.
Web: \href{http://www.stat.washington.edu/hoff/}{\tt http://www.stat.washington.edu/hoff/}. The author thanks Thomas Richardson for helpful discussions.}
        }
\date{ \today }
\maketitle

\begin{abstract}
Quantitative studies in many fields involve the analysis 
of multivariate data of diverse types,
 including measurements that we may consider 
binary, ordinal and continuous.  
One approach to the analysis of such mixed data
is to use a copula model, in which  the associations
among the variables are parameterized separately from
 their univariate marginal distributions.
The purpose of this article is to provide a simple, general method of semiparametric inference
for copula models
via a type of rank likelihood function for the
association parameters.
The proposed method of inference can be viewed as a generalization of
marginal likelihood estimation, in which
inference for a parameter of interest is based on a summary statistic
whose sampling distribution is not a function of any nuisance parameters.
In the context of copula estimation, the extended rank likelihood
is a function of the association parameters
only and its applicability does not depend on any assumptions about the marginal
 distributions of the data,
thus making it
appropriate for the
analysis of mixed continuous and discrete data with arbitrary
marginal distributions.
Estimation and inference for parameters of the Gaussian copula are available
via a straightforward
Markov chain Monte Carlo algorithm based on Gibbs sampling.
Specification of prior distributions  or a parametric form for the univariate marginal
distributions of the data is not necessary.

\vspace{.2in}
\noindent {\it Some key words}:
Bayesian inference, latent variable model,
marginal likelihood, Markov chain Monte Carlo,
multivariate estimation, polychoric correlation,
rank likelihood, sufficiency.
\end{abstract}

\section{Introduction}
Studies involving multivariate data often include measurements of
diverse types. 
For example, a survey or observational study may 
record the sex,  
education level and income of its participants, 
thus including measurements that we may consider
binary, ordinal and continuous.  
Such studies are generally concerned with  statistical 
associations among the variables, but not necessarily the scale on which
the variables  are measured. One approach to data analysis  in these situations is to obtain 
rank-based measures of bivariate association, such as the 
rank correlation or ``Spearman's rho''.  Such procedures
are scale-free, but involve ad-hoc methods for dealing with ties
 and provide inference that is generally limited
to hypothesis tests of bivariate association.
These issues make such procedures problematic for the analysis of 
much of social science survey data, in which the 
variables are often discrete  
and the hypotheses of interest generally concern multivariate and 
conditional associations. 
\begin{figure}
\centerline{\includegraphics[height=4.5in]{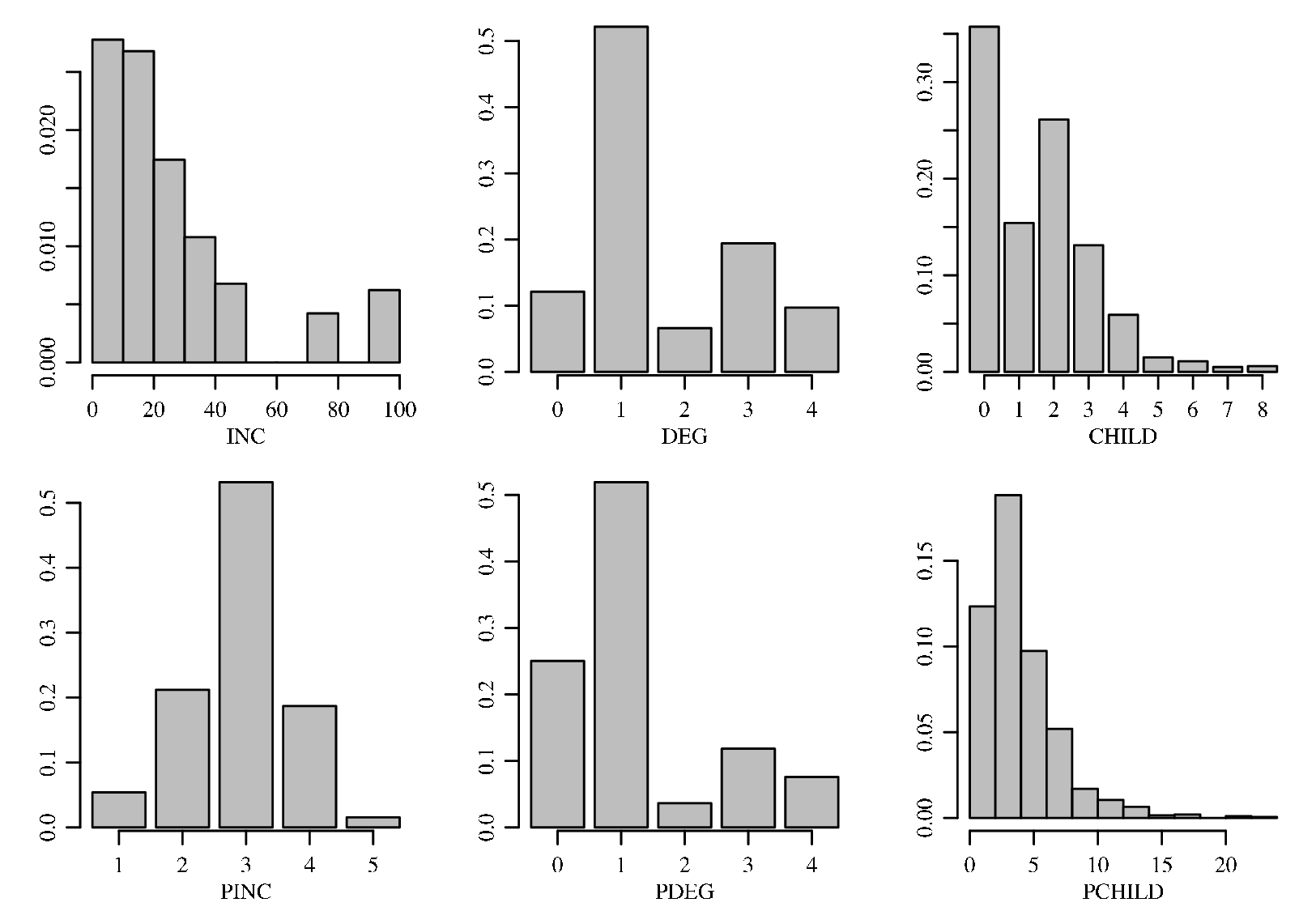}}
\caption{Univariate histograms of the GSS data.}
\label{fig:udist}
\end{figure}
 For example, 
Figure \ref{fig:udist} shows histograms of six demographic
variables of male respondents to the 1994 General Social Survey.  
The variables INC, DEG and CHILD refer to 
the income, highest degree and number of children
of  a survey respondent, 
and
PINC, PDEG and PCHILD refer to similar variables
of the respondent's parents (further details on the variables are
given in Section 4). 
All of these variables are ordered categorical variables, even 
though some of them have many levels. 
Additionally, our interests in these variables 
 involve measures of conditional association: 
An assessment of the relationship
between 
income  and number of children would generally be considered 
incomplete if it failed to 
account for  heterogeneity of the survey respondents 
 in terms of their age, parental  income
and other variables. 

The standard approach to making conditional assessments of statistical 
association is the use of regression models. For example, to describe 
the conditional association between income and number of children we could 
estimate  the parameters in a   regression  model of the following form:
\begin{equation}  {\rm INC}_i = \beta_0 + 
   \beta_1 {\rm CHILD}_i +\beta_2 {\rm DEG}_i + \beta_3 
   {\rm AGE}_i +  
 \beta_4 {\rm PCHILD}_i + \beta_5 {\rm PINC}_i +\beta_6 {\rm PDEG}_i + 
 \epsilon_i 
\label{eq:conmod1}
\end{equation}
Least-squares parameter estimates for this model, along with normal-theory
 $p$-values  appear in the first row of Table \ref{tab:conmod}. 
Standard practice is to interpret the   
$p$-value 
of $0.11$ for CHILD as suggesting  that there is not substantial evidence 
against $\beta_1=0$, in which case the 
model implies that
INC and CHILD are conditionally 
independent given the other variables. Alternatively, 
we could have evaluated the same conditional independence hypothesis 
with a regression model for CHILD. As this is a count variable, we might 
 use 
a Poisson regression model:
\begin{equation}  {\rm CHILD}_i \sim  {\rm Pois}( \exp\{ \beta_0 + 
   \beta_1 {\rm INC}_i +\beta_2 {\rm DEG}_i + \beta_3 
   {\rm AGE}_i +  
 \beta_4 {\rm PCHILD}_i + \beta_5 {\rm PINC}_i +\beta_6 {\rm PDEG}_i \})
\label{eq:conmod2}
\end{equation}
Maximum likelihood estimates and $p$-values for this model appear in 
the second row of Table \ref{tab:conmod}. In contrast to the results of Model 
(\ref{eq:conmod1}), these results indicate reasonably strong evidence ($p=0.01$) 
that CHILD and INC are \emph{not} conditionally independent, given the other variables. 

\begin{table}
\begin{small}
\begin{center}
\begin{tabular}{r||r|r|r|r|r|r|r}
&  \multicolumn{7}{c}{Predictor} \\
Response & INC & CHILD & DEG & AGE & PCHILD & PINC & PDEG  \\ \hline \hline
INC & NA & {\bf  1.10 (.11)} &   7.03 ($<$.01)&   .34 ($<$.01)&   4.07
( $<$.01) &  .28   (.41 ) &  1.40  ( .12) \\ \hline
CHILD &  {\bf  .01 (.01)} & NA & -.07  (.06) &  .04 ($<$.01)& -.06 (.20)&
  .02  (.08) & -.05  (.20)
\end{tabular}
\end{center}
\label{tab:conmod}
\end{small}
\caption{Estimation linear and Poisson regression coefficients  in the 
conditional models for INC and CHILD, 
with $p$-values   in parentheses.  }
\end{table}

The contradiction between the above two analyses 
is partly 
due to the inadequacies of the simple univariate parametric   Gaussian
and Poisson models. However, in general there is no reason to expect 
that two separately estimated conditional 
models will give compatible results: Given two conditional models
$f_1(y_1|y_2,\m x)$ and $f_2(y_2|y_1,\m x)$, only under very specific 
conditions does there exist a joint 
probability 
distribution $p(y_1,y_2|\m x) $ having  $f_1$ and $f_2$ as its  full conditional
distributions \citep{arnold_press_1989}. 
This presents a problem for the analysis of multivariate data 
of diverse types:
in the absence  of an appropriate multivariate model, 
 common practice is to analyze the data via
one or more
univariate regression models, choosing the ``response'' from the variables
which might best fit an ordinary or generalized linear
regression model. However, as the above example shows, 
different choices about which variables to treat as the response can 
lead to incompatible models with different conclusions.

Part of the above problem can be resolved by jointly modeling 
the variables of interest. 
A number of latent-variable methods have been recently 
developed to accommodate non-Gaussian multivariate data. 
These methods 
generally
proceed by modeling
each component of a vector of observations with a parametric
exponential family model, in which the parameters  for each component
involve an unobserved latent variable.
For example,
\citet{chib_winkelmann_2001}
present a model for a vector of
correlated count data in which each component
is a Poisson random variable with a mean depending on a component-specific
latent variable. Dependence among the count variables is induced
by modeling the vector of latent variables with a multivariate
normal distribution. Similar approaches are proposed by
\citet{dunson_2000}
and described in Chapter 8 of \citet{congdon_2003}. 
The model of Chib and Winkelmann can be viewed as a  copula model,
in which the association parameters are modeled 
separately from the marginal distributions of the observed data.
Such  a modeling  approach  can be applied to  a wide variety
of multivariate analysis problems:
An old mathematical result known as Sklar's Theorem says
that every multivariate
probability distribution can be represented by its univariate
marginal distributions and a copula, which is a type of joint
distribution with fixed marginals.

\citet{pitt_chan_kohn_2006} develop an estimation procedure 
for multivariate normal copula models  in which the marginal distributions
belong to specified parametric families. 
Unfortunately, 
the marginal distributions of 
survey data such as age, number of children, income and
education level generally do not belong to standard
families.
For such data a 
semiparametric estimation strategy may be appropriate, in which
the associations among the variables are represented with
a simple parametric model but the marginal distributions
are estimated nonparametrically. 
In the case where all the variables are continuous,
\citet{genest_ghoudi_rivest_1995}
suggest a ``pseudo-likelihood'' approach to estimation, 
in which the observed data is transformed via the empirical marginal
distributions to obtain pseudo-data that can be used to estimate
the association parameters. 
\citet{klassen_wellner_1997}
study a similar type of 
estimation in the case of the Gaussian copula.  
Such estimators are well-behaved for continuous data but can fail for discrete
data, making them somewhat inappropriate for the analysis of
mixed continuous and discrete data. 
For ordinal discrete data with a known number of categories, the
dependence induced by the Gaussian copula model is called
polychoric correlation.
\citet{olsson_1979} describes a two-stage
estimation procedure for the parameters in the copula, and
this and other estimation strategies appear in a number
of software packages including SAS PROC FREQ  and the LISREL
module PRELIS.
\citet{kottas_muller_quintana_05} describe a nonparametric
estimation procedure in which the copula is based on a
mixture of normal distributions.
However, such procedures do not accommodate continuous data,
and may even be problematic for
discrete data with a large number of categories, as
inference in this case requires the simultaneous estimation of the
large number of parameters
specifying the marginal distributions.

As an alternative to these procedures, this article presents an approach
to copula estimation in which the marginal distributions are
arbitrary and of unspecified types, thus accommodating
both discrete and continuous data.
This is achieved by the use of a likelihood function
that depends on the association parameters
only, and does not make assumptions about the form of the
univariate marginal distributions.
Inference based on such a likelihood
is therefore
appropriate for the
joint analysis of continuous and ordinal discrete data.
For continuous data, the likelihood function we propose is
derived from the marginal probability of the ranks,
and can be seen as a multivariate version of a ``rank likelihood''
\citep{pettitt_1982,heller_qin_2001} which does not depend on
the univariate marginal distributions.
Unfortunately,
for discrete data the probability of the observed ranks is not free of
these nuisance parameters.
To solve this problem,
we derive a likelihood that is equivalent to the
distribution of the ranks for continuous data but is also
free of the nuisance parameters for discrete data.
This likelihood function is derived from  the probability
that the latent variables of the copula model satisfy the  partial
ordering induced by the observed data.
We call this function an extended rank likelihood, as it
generalizes the concept of rank likelihood. 
This likelihood can also  be seen as a generalization of a marginal likelihood,
which  is based on
a statistic whose sampling distribution
depends only on the parameter of interest and not on any nuisance parameters.



In what follows we work with the Gaussian copula model, although
the basic ideas can be extended to other parametric families of copulas.
In the next section we review  the general Gaussian copula model, and
discuss how inference for discrete data using existing semiparametric
methods is problematic.
Section 3 derives the extended rank likelihood as a general approach
to semiparametric copula estimation
and  discusses parameter estimation in the context of Bayesian
inference using a relatively simple Gibbs sampling scheme.

The primary goal of this paper is to provide 
a simple method of inference for the multivariate 
relationships between variables, such as INC, CHILD, DEG 
described above, whose univariate marginal distributions 
cannot  be well approximated with simple parametric 
models. 
In Section 4 we present an analysis of 
these and other demographic characteristics of males in the 1994 U.S.\ 
workforce and 
their parents. 
In particular, we are interested in 
the statistical associations among income, education and 
number of children of the survey respondents, and how they relate to 
similar characteristics of the parents of the survey 
respondents. 
The data come from the 1994 General Social Survey, and include
 a number of
discrete and non-Gaussian random variables. 
In addition to estimating a Gaussian copula model for these data, 
we estimate and describe the conditional dependencies  among the variables on the Gaussian 
scale, as well as provide predictive  and conditional 
distributions on the original 
scale  of the data.

Section 5 considers notions of statistical sufficiency relevant to the
 rank likelihood, and a
discussion follows in Section 6.

\section{Semiparametric copula estimation}
Let $y_1$ and $y_2$ be two random variables with continuous CDF's
$F_1$ and $F_2$. The transformed variables
$u_1 = F_1(y_1)$ and $u_2=F_2(y_2)$ both have uniform
marginal distributions.
The term ``copula modeling'' generally refers to a model
that parametrizes the joint distribution of
$u_1$ and $u_2$ separately from the marginal distributions
$F_1$ and $F_2$.
A semiparametric copula model
includes
a parametric model for the joint distribution of
$u_1$ and $u_2$, but lacks  
any parametric restrictions on $F_1$ or $F_2$.

Any continuous multivariate distribution
 can be used to form  
a copula model via an inverse-CDF transformation.  
For example, the bivariate normal distribution can be used to
generate dependent data with arbitrary marginals
$F_1$ and $F_2$ as follows:
\begin{enumerate}
\item sample $\left (\begin{array}{c}z_1 \\ z_2 \end{array} \right )  \sim $ 
bivariate normal$\left [ \left (\begin{array}{c}0 \\ 0 \end{array} \right ), 
 \left ( \begin{array}{cc} 1 & \rho \\ \rho & 1 \end{array} \right ) \right ]$;
\item set $y_1 = F_1^{-1}[\Phi(z_1)]$,  $y_2 = F_2^{-1}[\Phi(z_2)]$,
\end{enumerate}
where
$F^{-1}(u) =\inf  \{  y: F(y)\geq u\}$ denotes the
pseudo-inverse of a CDF $F$. 
The correspondence to the usual copula formulation can be seen
by noting that $\Phi(z) = u$ is uniformly distributed.

Suppose  $(y_{1,1},y_{1,2}),\ldots, 
 (y_{n,1},y_{n,2})$  are samples from a population that we wish to
model with a Gaussian copula. 
If the marginal distributions $F_1$ and $F_2$ were continuous and known,
then the values  $z_{i,j} = \Phi^{-1}[F_j(y_{i,j})]$ could be treated
as observed data and $\rho$ could be estimated directly from
the $z$'s, perhaps using  
the  unbiased estimator
$\hat \rho = \frac{1}{n}\sum_{i=1}^n  z_{i,1}z_{i,2} $.
Of course, the marginal CDF's are not typically known.
One semiparametric estimation strategy
is to plug-in the 
the empirical CDF's $\hat F_1$ and $\hat F_2$ to obtain 
pseudo-data $\tilde z_{i,j}= \Phi^{-1}[\frac{n}{n+1}\hat F_j(y_{i,j})] \equiv
   \Phi^{-1}[ \tilde F_{j}(y_{i,j})]$, where the rescaling is to avoid
infinities. 
For continuous data, the estimator $\tilde \rho = \frac{1}{n}\sum_{i=1}^n  \tilde z_{i,1}\tilde z_{i,2} $
is asymptotically equivalent to the asymptotically efficient Van der Waerden
normal-scores rank correlation coefficient 
\citep{hajek_sidek_1967,klassen_wellner_1997}. 
This estimator is similar to one obtained 
from a more general  pseudo-likelihood estimation procedure
described and studied by
\citet{genest_ghoudi_rivest_1995}. 
In the context of the Gaussian copula model, 
the maximum pseudo-likelihood
procedure is to 
\begin{enumerate}
\item set $\tilde z_{i,j} = \Phi^{-1}[\tilde F_j(y_{i,j})];$ 
\item maximize in $\rho$ the  pseudo-log-likelihood  $\sum_{i=1}^n 
       \log {\rm bvn}(\tilde z_{i,1},\tilde z_{i,2}|\rho)$,
\end{enumerate}
where ${\rm bvn}(\cdot |\rho)$ denotes the bivariate normal density with
standard normal marginals. 
Genest et al.\  show that the resulting pseudo-likelihood estimator
is consistent and asymptotically normal under the condition
that $F_1$ and $F_2$ are continuous. 
However, 
this condition calls into question the appropriateness
of the pseudo-likelihood approach for non-continuous data 
such as sex, education level, age or any other type of data
where there are likely to be ties.

What could go wrong with such an estimator in situations involving
discrete data?
In general, these pseudo-data estimators of copula parameters 
will be problematic for discrete data because transformations of
such  data do not really change the data distribution, they
just change the sample space. 
Consider the simple case of a continuous variable $y_1$ and 
a binary variable $y_2$ such that $\Pr(y_2=0) =\Pr(y_2=1) = 1/2$.
Letting $\tilde z_{i,j} = \Phi^{-1}[\tilde F_j(y_{i,j})]$,  
the distribution of $\tilde z_{1,1},\ldots, \tilde z_{n,1}$ will have
an approximately standard normal distribution, but
$\tilde z_{i,2}$ will be approximately equal 
to either  $\Phi^{-1}(\frac{1}{2} \frac{n}{n+1})$ or 
 $\Phi^{-1}(\frac{n}{n+1})$  with probability one-half each.
If the Gaussian copula model is correct, then 
one can show that the expectation of $\tilde \rho$
is roughly 
 $\frac{\rho}{\sqrt{2\pi}} \Phi^{-1}(\frac{n}{n+1})$.
As $n$ increases so does
the expectation of $\tilde \rho$, and it is not a consistent estimator.
One problem here is that all of the $\tilde z_{i,2}$'s such that
$y_{i,2}=1$ are being pushed to the extreme standard normal
quantile $\Phi^{-1}(\frac{n}{n+1})$, which 
in the case of continuous data
would happen just to  a single datapoint.
The situation is only partly improved by using the sample
correlation of the pseudo-data as an estimator:
The variance of $\tilde z_1$ is
approximately 1 and the variance of $\tilde z_2$ is approximately
 $[ \frac{1}{2}\Phi^{-1} (\frac{n}{n+1}) ]^2$, giving an approximate
sample correlation of ${\rm Cor}(\tilde z_{i,1},\tilde z_{i,2}) \approx
 \rho \sqrt{2/\pi}$.

\section{Estimation using the extended rank likelihood}
In this section we derive a likelihood function that depends on
the association parameters and not on the unknown marginal distributions.
For continuous data this function is equivalent to 
the distribution of the multivariate ranks.
This is not the case of discrete data, for which
the distribution of the  ranks
 depends on the univariate marginal distributions.
In this case  the derived likelihood
function contains less total information than one based on
the  ranks,
 but it is free of any parameters describing the
marginal distributions.



\subsection{Extended rank likelihood}
Generalizing from the previous section,
the Gaussian copula sampling model can
be expressed as follows:
\begin{eqnarray}
\m z_1,\ldots,\m z_n|\m C&\sim&\mbox{i.i.d.\ multivariate normal}(\m 0, \m C), 
   \\  \label{eq:gcm}
y_{i,j} &=& F^{-1}_j [ \Phi(z_{i,j}) ], \nonumber
\end{eqnarray}
where $\m C$ is a $p\times p$ correlation matrix and
each $F_j^{-1}$ denotes the (pseudo) inverse of an unknown
univariate
CDF,
not necessarily
continuous.

Our goal is to make inference on $\m C$, and not on
the potentially high-dimensional parameters $F_1,\ldots, F_p$.
If the $\m z$'s were observed  we could use them
to directly estimate $\m C$. The
$\m z$'s are not observed of course, but the
$\m y$'s do provide a limited amount of information about
them,
even absent any knowledge of
 the $F$'s:
Since the $F$'s are non-decreasing,
observing  $y_{i_1,j}<y_{i_2,j}$ implies that
$z_{i_1,j} < z_{i_2,j}$.
More generally, observing $\m Y=( \m y_1,\ldots, \m y_n)^T$ tells us that
$\m Z = ( \m z_1,\ldots, \m z_n)^T$
must lie in the set
\[ 
\{ \m Z \in \mathbb R^{n\times p} :  
\max\{z_{k,j}:y_{k,j}<y_{i,j} \} < z_{i,j} <   \min\{z_{k,j}:y_{i,j}<y_{k,j}\} \}. \]
We can take the occurrence of this event as our data.
Letting $D$ be the fixed subset  of $\mathbb R^{n\times p}$
generated by the observed value of $\m Y$,
we can calculate the following ``likelihood'':
\begin{equation} 
\Pr ( \m Z \in   D | \m C , F_1,\ldots, F_p) 
  = \int_{D} p(\m Z | \m C) \ d\m Z 
               = \Pr ( \m Z \in   D| \m C ) . 
\label{eq:msl}
\end{equation}
As a function of the parameters, this likelihood depends only
on the  parameter of interest $\m C$ and not the
 nuisance parameters
$F_1,\ldots, F_p$.
Estimation of $\m C$ can proceed by maximizing $\Pr(\m Z\in D|\m C) $
as a function of $\m C$, or by obtaining a posterior distribution
$\Pr(\m C | \m Z \in D) \propto p(\m C) \times \Pr(\m Z\in D | \m C)$.

The likelihood function (\ref{eq:msl}) can be seen as a type of marginal likelihood
function for estimation in the presence of a nuisance parameter:
Consider a generic statistical problem in which the density for data $y$
depends on a parameter of interest $\theta$ and a nuisance parameter
$\psi$. If there exists a statistic $t(y)$ whose distribution depends
on $\theta$ only, then the density of $y$ may be decomposed as
\begin{eqnarray*}
 p(y|\theta,\psi) &=& p(t(y),y | \theta,\psi) \\
 &=& p(t(y)|\theta) \times p(y|t(y),\theta,\psi). 
\end{eqnarray*}
In this situation, estimation of $\theta$ can be based on
the marginal likelihood $p(t(y)|\theta)$,
eliminating the need to estimate the nuisance parameter $\psi$
(see, for example, Section 8.3 of \citet{severini_2000}). 
The likelihood function  $\Pr(\m Z\in D|\m C)$ in our copula estimation
problem can be derived  analogously, by
decomposing  the probability
of the observed data as
\begin{eqnarray}
p(\m Y|\m C,F_1,\ldots,F_p) &=& p(\m Z \in D, \m Y  |\m C, F_1,\ldots,F_p) 
 \label{eq:pr1} \\
   &=& \Pr(\m Z \in D | \m C ) \times p( \m Y | \m Z \in D, \m C, 
              F_1,\ldots,F_p ) . \label{eq:pr2}
\end{eqnarray}
Equation (\ref{eq:pr1}) holds because
the event
 $\m Z \in D$  occurs whenever $\m Y$ is observed.
This derivation can be made rigorous by
deriving the density $p(\m Y|\m C, F_1,\ldots, F_p)$ from the limit
of $\Pr( \cap_{i,j} (y_{i,j}-\epsilon,y_{i,j}]|\m C, F_1,\ldots, F_p)$
as $\epsilon \rightarrow 0$.
As in the case of marginal likelihood, our approach is to estimate $\m C$
using only  $\Pr(\m Z\in D|\m C)$, the part of the
observed data likelihood (\ref{eq:pr2})
that depends on the parameter of interest $\m C$ and not on the nuisance
parameters $F_1,\ldots, F_p$. Since our likelihood function is
based on the marginal probability of an event that is a 
superset of observing the ranks,
we refer to it as an extended rank likelihood.


\subsection{Estimation of the copula parameters}
Bayesian inference for $\m C$
can be achieved via construction of a
Markov chain having a stationary distribution
equal to $p(\m C|\m Z\in D) \propto p(\m C) \times p(\m Z \in D | \m C)$.
In the case of the Gaussian copula with a semi-conjugate prior distribution,
the Markov chain can be constructed quite easily
using
Gibbs sampling.
This prior distribution
for $\m C$
is defined as follows:
Let $\m V$ have an inverse-Wishart$( \nu_0, \nu_0 \m V_0)$ prior distribution,
parameterized so that $E[\m V^{-1}] = \m V_0^{-1}$,
and let $\m C$ be  equal in distribution to the
the correlation matrix with entries
$\m V_{[i,j]}/\sqrt{ \m V_{[i,i]} \m V_{[j,j]} }$.
Using this prior distribution,
approximate samples from $p(\m C|\m Z\in D)$ can be obtained
by iterating the following Gibbs sampling scheme:

\begin{description}
\item[Resample $\m Z$.] Iteratively over $(i,j)$, sample
     $z_{i,j}$ from $p(z_{i,j}|\m V, \m Z_{[-i,-j]}, \m Z\in D)$ as follows:

\begin{itemize}
\item[] For each  $j \in \{1,\ldots, p\}$
\begin{itemize}
\item[] For each $y \in$ unique$\{ y_{1,j},\ldots, y_{n,j}\}$
\begin{itemize}
\item[1.] Compute $z_l = \max \{ z_{i,j} : y_{i,j} < y \}$ and
              $z_u = \min \{ z_{i,j} :   y< y_{i,j}  \}$
\item[2.] For each $i$ such that $y_{i,j}=y$,
\begin{itemize}
\item[(a)] compute $\sigma_j^2  =  \m V_{[j,j]} - \m V_{[j,-j]}  \m V_{[-j,-j]}^{-1} \m V_{[-j,j]} $
\item[(b)]  compute $\mu_{i,j} = \m Z_{[i,-j]} ( \m V_{[j,-j]} \m V_{[-j,-j]}^{-1} )^T$
\item[(c)] Sample $u_{i,j}$ uniformly from $(\Phi[\frac{z_l-\mu_{i,j}}{\sigma_j} ],
   \Phi[\frac{z_u-\mu_{i,j}}{\sigma_j} ] )$
\item[(d)] Set $z_{i,j} = \mu_{i,j} + \sigma_j \times \Phi^{-1}(u_{i,j}) $
\end{itemize}
\end{itemize}
\end{itemize}
\end{itemize}

\item[Resample $\m V$.] Sample $\m V$ from an inverse-Wishart$(\nu_0+n, 
       \nu_0\m V_0 + \m Z^T \m Z )$ distribution.
\item[Compute $\m C$.] Let $\m C_{[i,j]} =  \m V_{[i,j]}/\sqrt{ \m V_{[i,i]} \m V_{[j,j]} }$.
\end{description}
Iteration of this algorithm generates a Markov chain in $\m C$ whose
stationary distribution is $p(\m C| \m Z \in D)$.
This algorithm is easily modified to accommodate data that are
missing-at-random: If $y_{i,j}$ is missing,
the full conditional distribution of $z_{i,j}$
is the unconstrained normal distribution with mean $\mu_{i,j}$ and
variance $\sigma^2_j$ given above.

The 
reader may have noticed that the samples of $\m Z$ are
based on the covariance matrix $\m V$ and not the correlation matrix $\m C$.
To see why this does not matter for estimation of $\m C$, compare our
original model,
\begin{eqnarray*}
\m V & \sim &  \mbox{inverse-Wishart}(\nu_0, \nu_0 \m V_0)  \\ 
\{ \m C_{[i,j]}  \} &=& \{ \m V_{[i,j]}/\sqrt{ \m V_{[i,i]} \m V_{[j,j]} } \} \\
\m z_1,\ldots, \m z_n &\sim& \mbox{i.i.d.\ multivariate normal}(\m 0, \m C) \\
y_{i,j} &=& G_j (z_{i,j} ), 
\end{eqnarray*}
to the equivalent model
\begin{eqnarray*}
\m V & \sim &  \mbox{inverse-Wishart}(\nu_0, \nu_0 \m V_0)  \\ 
\m z_1,\ldots, \m z_n &\sim& \mbox{i.i.d.\ multivariate normal}(\m 0, \m V) \\
\tilde z_{i,j} &=& z_{i,j}/\sqrt{\m V_{[j,j]}}, \ \mbox{ and let } \m C = {\rm Cov}(
 \tilde {\m z} )     \\
y_{i,j} &=& G_j (\tilde z_{i,j} ). 
\end{eqnarray*}
The  $\m z$'s  in the first formulation are equal in distribution
to the $\tilde {\m z} $'s  in the second,
and so posterior inference for $\m C$  is
equivalent under either model.
The Gibbs sampling scheme outlined above is based on a Markov chain
in $\m V$ and $\m z_1,\ldots, \m z_n$  based on the second formulation.
Note that in this formulation the observed data implies the same ordering
$D$ on both the $\tilde {\m z}$'s and the $\m z$'s.
Additionally, posterior estimation of $\m C$ is invariant to 
changes in the prior distribution on $\m V$ that do not alter 
the induced prior on $\m C$. For example, if $\m V_0$ and $\m V_0'$
are two different covariance matrices with the same correlations, 
then the posterior distribution of $\m C$ under $\m V \sim$ inverse-Wishart$(\nu_0, \nu_0 \m V_0)$ will be equal to that under 
$\m V \sim$ inverse-Wishart$(\nu_0, \nu_0 \m V_0')$ .



\section{Income, education and intergenerational  mobility}
The U.S.\ census reports  a strong positive relationship between income 
and educational attainment 
\citep{day_newberger_2002}. 
However, in many studies both of these variables have been shown to be
associated with a number of family background variables such as parental income,
parental educational attainment and number of siblings
\citep{ermisch_francesconi_2001, blake_1985}. 
Additionally, 
some researchers have suggested that having children reduces 
opportunities for educational attainment
\citep{moore_waite_1977}, while others 
have found evidence that  economic status of males is positively associated
with their  fertility \citep{hopcroft_2006}. 
Results such as these are generally based on univariate 
regression models in which one variable from 
a sample survey is selected as a ``response'' or ``dependent'' variable 
and the 
others as ``control'' or ``independent'' variables. However, all 
of the variables in these studies are randomly sampled and all are 
potentially dependent 
on one another.

In this section we describe the multivariate dependencies among 
income, education and number of children using the 
Gaussian copula  model and the semiparametric estimation procedure
described in Section 3. 
Specifically, we analyze survey data on 1002 males in the U.S.\  labor force 
(meaning not retired, in school or in an institution), 
obtained from the 1994 General Social Survey. 
Data and details for the survey are available at
\href{http://webapp.icpsr.umich.edu/GSS/}{\tt http://webapp.icpsr.umich.edu/GSS/}.

The relevant variables for this analysis include
the income, education, and number of children of the survey respondent, 
as well as similar variables for the respondent's parents. Age 
of the survey respondent is additionally
included, as it  is typically 
strongly related to income and number of children. 
The measurement scales for these variables are as follows:

\medskip

\begin{tabular}{rl}
INC:  & income of the respondent in 1000s of dollars, binned into 21 
ordered     categories \\
DEG: & highest degree ever obtained (None, HS, Associates, Bachelors, Graduate)\\
CHILD:  &  number of children ever had \\
PINC: & financial status of respondent's parents when respondent was 16
 (on a 5-point scale) \\
PDEG: &  maximum of mother's and father's highest degree\\
PCHILD:  & number of siblings of the respondent plus one \\
AGE: & age of the respondent in years \\
\end{tabular}

\medskip
Missing data rates among each of the non-income variables was less
than 4\%. The missing data rates for INC and PINC were  10\% and  
 48\% respectively.  However, the question PINC was asked on 
only half of the surveys, and so missing values for this variable can
 reasonably
be considered as missing at random.

\begin{figure}
\centerline{\includegraphics[height=3.5in]{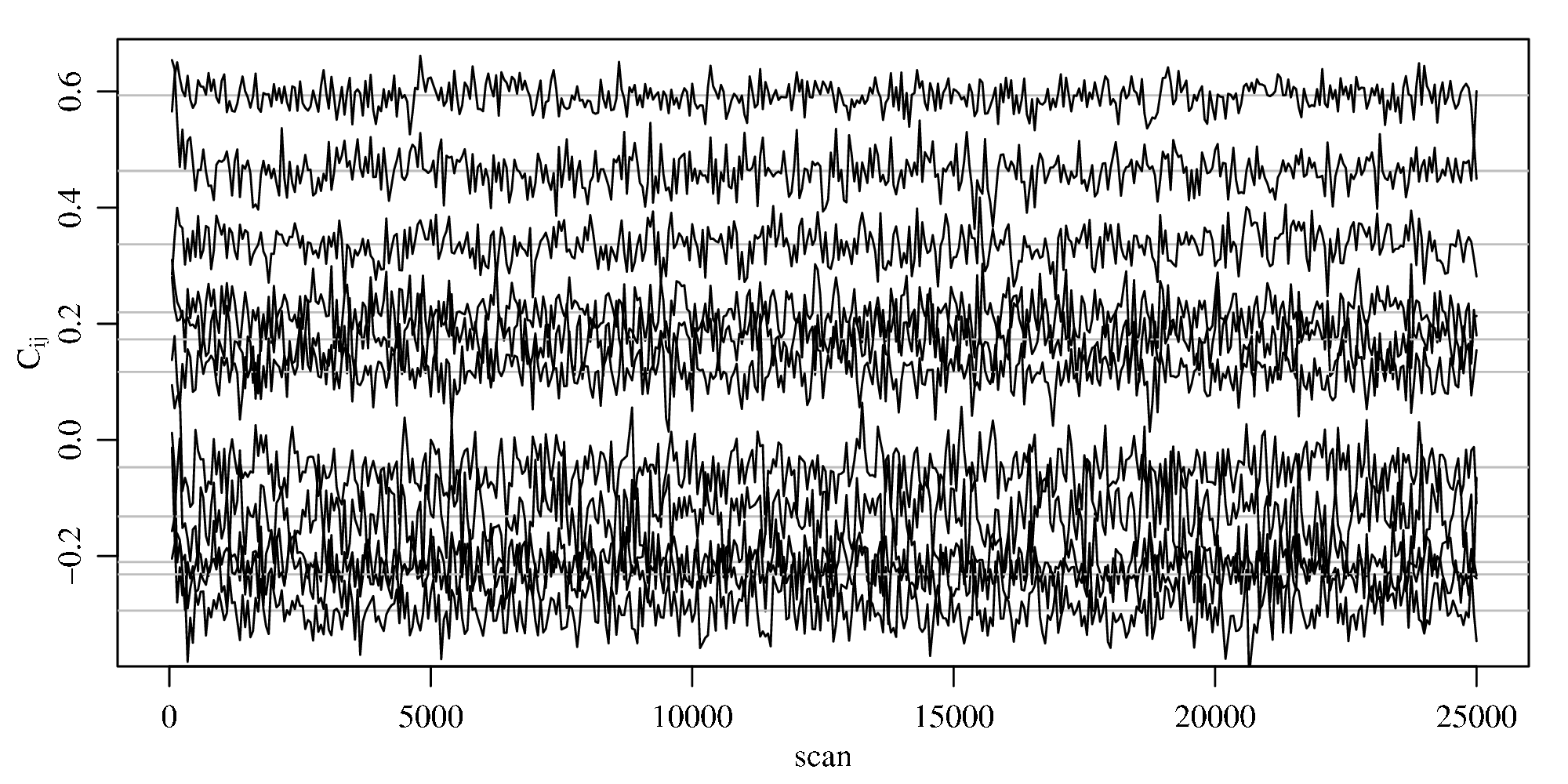}}
\caption{MCMC samples of 11 of the
correlation coefficients, plotted every 50th scan.}
\label{fig:mcmcdiag}
\end{figure}

\subsection{Estimation of $\m C$}
Using an inverse-Wishart $( p+2, (p+2)\times \m I)$ prior distribution
for $\m V$, the Gibbs sampling scheme outlined in Section 3
was iterated 25,000 times with parameter values saved
 every 10 scans, resulting in
2500 samples of $\m C$ for posterior analysis. Mixing of the Markov chain
was quite good: Figure \ref{fig:mcmcdiag} shows MCMC samples of 11 
elements  of $\m C$, 
 corresponding to the odd order statistics of 
$E[\m C|\m Z \in D]$. 
Convergence to stationarity appears to occur quickly, almost certainly within
the first 5000 scans. 
Dropping these scans to allow for burn-in,
we are left with 2000 saved scans for posterior analysis.
The autocorrelation across these saved scans was low, with the
lag-10 autocorrelation less
than 0.05 in absolute value for all elements of $\m C$, and
much closer to zero for most. Based on the autocorrelation in the Markov chain, 
the effective sample sizes for estimating the posterior means of the 
elements of $\m C$ were at least 1500.

\subsection{Posterior inference}
\begin{figure}
\centerline{\includegraphics[height=3.25in]{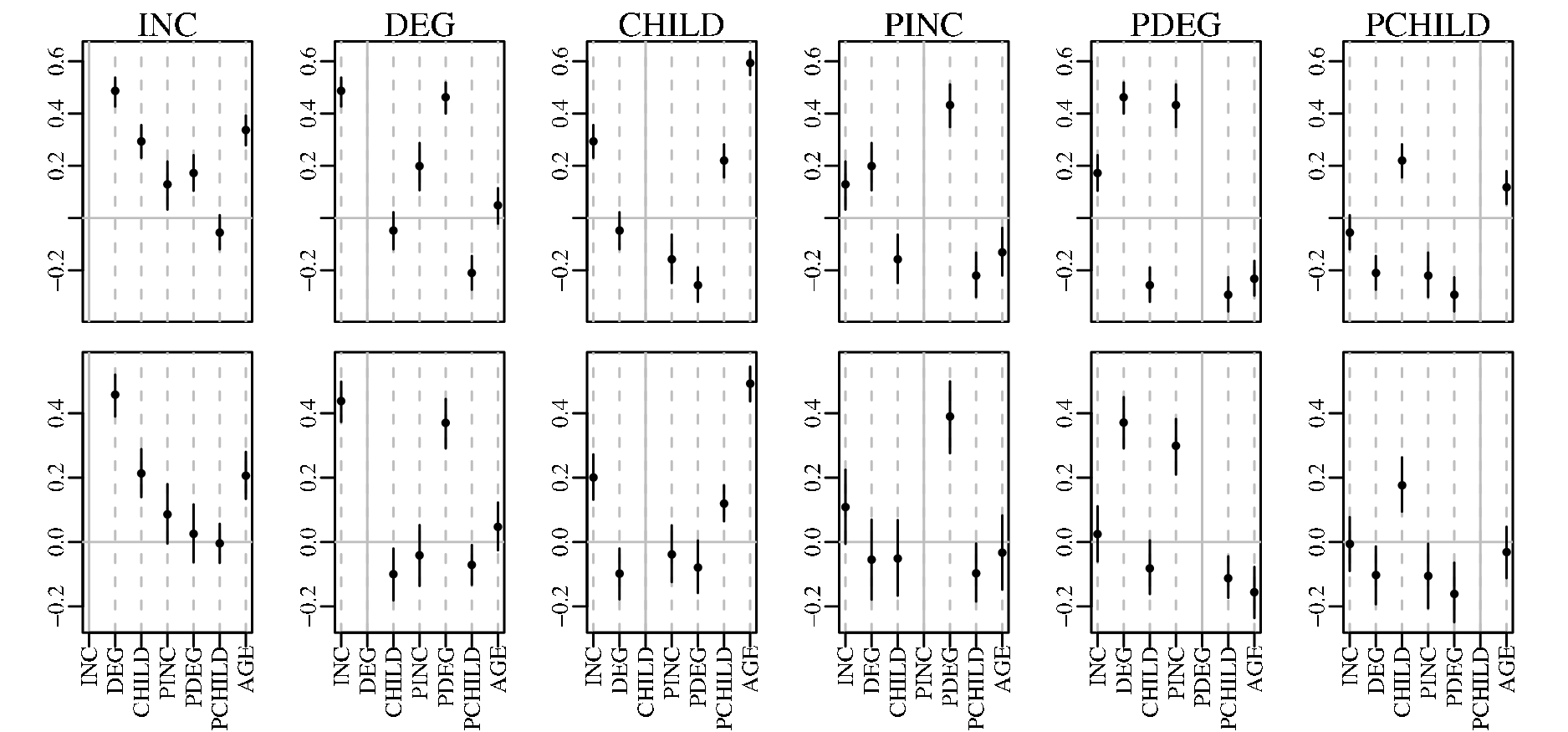}}
\caption{Dependence parameters for the
GSS data. The first row gives 2.5\%, 50\%, 97.5\% posterior
quantiles of the
correlation coefficients $E[z_j z_k]$.
The second row  gives the regression coefficients
  $\nabla E[ z_j |z_{-j}]$. }
\label{fig:postmean}
\end{figure}
Posterior distributions of the correlation parameters
are summarized in the first and second rows of Figure  \ref{fig:postmean}. 
The first row  gives
 2.5\%, 50\% and 97.5\% posterior quantiles of
the correlation coefficients, representing 
scale-invariant bivariate associations among the six variables of interest. 
The fact that most of these 95\% credible intervals do not 
contain zero 
indicates that most variables are associated 
with most of the other variables. For example, the results 
suggest that
INC has non-zero positive correlations with 
DEG, CHILD, PINC, PDEG and AGE, and a weak negative 
correlation with PCHILD. 
DEG shows positive correlations with INC, PINC, PDEG, 
and negative correlation with PCHILD (in accordance with the 
conclusion  of \citet{blake_1985}).  

\begin{figure}
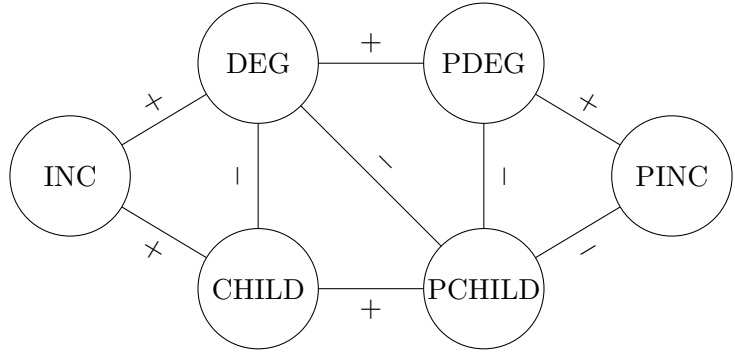

\begin{center}
\begin{pgfpicture}{-4.5cm}{-2.5cm}{4.5cm}{2.5cm} 
\pgfnodecircle{inc}[stroke]{\pgfxy(-4,0)}{.8cm}
\pgfnodecircle{chi}[stroke]{\pgfxy(-1.5,-1.5)}{.8cm} 
\pgfnodecircle{deg}[stroke]{\pgfxy(-1.5,1.5)}{.8cm} 
\pgfnodecircle{pdeg}[stroke]{\pgfxy(1.5,1.5)}{.8cm} 
\pgfnodecircle{pchi}[stroke]{\pgfxy(1.5,-1.5)}{.8cm} 
\pgfnodecircle{pinc}[stroke]{\pgfxy(4,0)}{.8cm} 

\pgfputat{\pgfnodecenter{inc}}{\pgfbox[center,center]{INC}}
\pgfputat{\pgfnodecenter{chi}}{\pgfbox[center,center]{CHILD}}
\pgfputat{\pgfnodecenter{deg}}{\pgfbox[center,center]{DEG}}
\pgfputat{\pgfnodecenter{pdeg}}{\pgfbox[center,center]{PDEG}}
\pgfputat{\pgfnodecenter{pchi}}{\pgfbox[center,center]{PCHILD}}
\pgfputat{\pgfnodecenter{pinc}}{\pgfbox[center,center]{PINC}}

\pgfnodeconnline{inc}{deg} 
\pgfnodeconnline{inc}{chi} 
\pgfnodeconnline{deg}{chi} 

\pgfnodeconnline{pinc}{pdeg}
\pgfnodeconnline{pinc}{pchi}
\pgfnodeconnline{pdeg}{pchi}

\pgfnodeconnline{deg}{pdeg}
\pgfnodeconnline{chi}{pchi}
\pgfnodeconnline{pchi}{deg}

\pgfnodelabelrotated{inc}{deg}[0.5][5pt]{\pgfbox[center,base]{$+$}} 
\pgfnodelabelrotated{chi}{inc}[0.5][5pt]{\pgfbox[center,base]{$+$}}     
\pgfnodelabelrotated{chi}{deg}[0.5][5pt]{\pgfbox[center,base]{$-$}}     
\pgfnodelabelrotated{deg}{pchi}[0.5][5pt]{\pgfbox[center,base]{$-$}}     
\pgfnodelabelrotated{deg}{pdeg}[0.5][5pt]{\pgfbox[center,base]{$+$}}     
\pgfnodelabelrotated{pchi}{chi}[0.5][5pt]{\pgfbox[center,base]{$+$}}     
\pgfnodelabelrotated{pdeg}{pinc}[0.5][5pt]{\pgfbox[center,base]{$+$}}     
\pgfnodelabelrotated{pinc}{pchi}[0.5][5pt]{\pgfbox[center,base]{$-$}}     
\pgfnodelabelrotated{pdeg}{pchi}[0.5][5pt]{\pgfbox[center,base]{$-$}}

\end{pgfpicture} 
\end{center}
\caption{Reduced conditional dependence graph for the GSS data.}
\label{fig:mdg}
\end{figure}

Perhaps of more interest are conditional associations. 
The second column of Figure \ref{fig:postmean}  gives the
2.5\%, 50\% and 97.5\% 
quantiles for the ``regression coefficients''
$ \m C_{[j,-j]} \m C_{[-j,-j]}^{-1}$ for each variable. 
These coefficients represent conditional dependencies among the 
underlying processes that give rise to the observed data. 
On this scale, the full conditional distribution of INC 
depends most strongly on DEG, and 
to a lesser extent on 
CHILD and AGE. Interestingly, the conditional relationship between
INC and PINC has a non-negligible ($>5\%$) probability of being  
less than or equal to zero.
Figure \ref{fig:mdg} summarizes these results with a graph 
indicating the conditional dependencies
among the $\m z$-variables corresponding to the  six variables of interest
(implicitly conditioning on AGE). An edge is present between 
two nodes if the 95\% credible interval for the associated 
regression parameter does not contain zero. 
This graph suggests that 
although INC and PINC are positively associated, this association 
is mediated by the intergenerational relationships of 
DEG, PDEG, CHILD and PCHILD.


\subsection{Conditional distributions for the 
INC, DEG, PINC relationship}
The results in Figure \ref{fig:postmean} suggest that, 
although INC and PINC are positively correlated, 
PINC is a relatively weak predictor of INC compared to DEG. 
However, PINC is a strong predictor of PDEG, and PDEG is a strong 
predictor of DEG, suggesting an indirect effect of PINC on INC. 

These  conclusions about INC, DEG and PINC 
are made in terms of associations among the $\m z$-variables, 
although  it is often desirable to report results on the 
scale of the original data. 
With this in mind, we now describe the relationship 
between INC, DEG and PINC on the original data scale, 
using  an estimated predictive distribution 
Pr(INC, DEG, PINC), which we decompose as $\Pr({\rm INC}|{\rm DEG},{\rm PINC}) \times \Pr({\rm DEG}|{\rm PINC}) \times \Pr({\rm PINC})$. 

A  predictive distribution for $\m y$ can be obtained in 
a few different ways. Perhaps the simplest method is to combine the 
posterior distribution of $\m C$ with the  empirical univariate
marginal distributions $\hat F_1,\ldots, \hat F_p$ of the observed 
data (an alternative method is presented in the Discussion). Using this method, a  predictive sample of $\m y$
 can be obtained as follows:
\begin{enumerate}
\item sample $\m C \sim p(\m C|\m Z\in D)$; 
\item sample $\m z \sim $ multivariate normal$(\m 0, \m C)$; 
\item set $y_j = \hat F^{-1}_j(z_j)$. 
\end{enumerate}
Although this somewhat ad-hoc approach disregards uncertainty in the estimation 
of $F_1,\ldots, F_p$ (for prediction of $\m y$, not for estimation of $\m C$), 
it provides a predictive joint 
distribution that matches the observed data in 
terms of the univariate marginal 
distributions but has a simple, smooth Gaussian copula representing 
multivariate dependence. From these predictive samples we can obtain 
Monte Carlo estimates of various quantities of interest, including 
a consistent set of 
conditional distributions on the original scale of the data. 

\begin{figure}
\centerline{\includegraphics[height=8.25in]{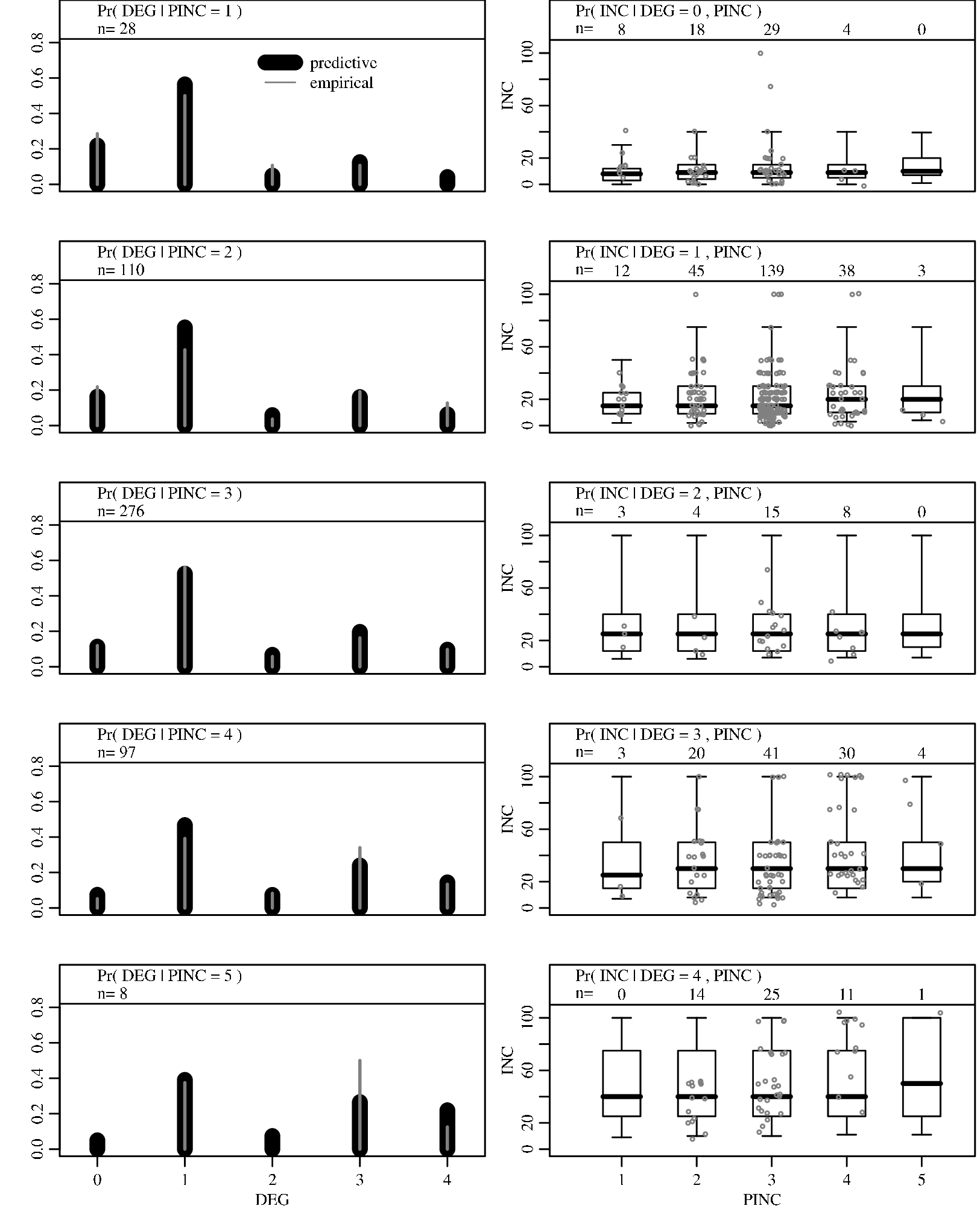}}
\caption{Empirical and  predictive conditional 
distributions for INC, DEG and PINC.}
\label{fig:cdist}
\end{figure}

The first column of Figure \ref{fig:cdist} plots the
predictive distribution of DEG conditional on ${\rm PINC}=x$  
for $x \in \{1,2,3,4,5\}$. 
As on the $\m z$-scale, large values of PINC 
correspond to large values of DEG. 
The estimated conditional probability of someone not finishing high-school
given PINC=5 is 5\%, whereas for PINC=1 it is 22\%, giving an 
odds ratio of odds(DEG=None$|$PINC=1) / odds(DEG=None$|$PINC=5) = 5.35. 
Similarly, the corresponding odds ratio for having a graduate degree is 
 odds(DEG=Grad$|$PINC=5) / odds(DEG=Grad$|$PINC=1) = 6.5.
For comparison, the empirical conditional distributions are 
provided on the same plot. In general the fit is good, 
with  most of the discrepancies occurring 
in categories of PINC with small sample sizes 
($n=28$ for PINC=1, and $n=8$ for PINC=5).
Note that if we were to estimate the above odds ratios using the empirical 
conditional distributions we would obtain ratios equal to infinity. 
In situations such as these where the sample size is low, we may prefer
to estimate conditional distributions with a model that can 
share information across the categories of  a variable, rather than use an 
empirical estimator having a high sampling variability. 

The second column of Figure {\ref{fig:cdist}} displays 
estimated  quantiles of $\Pr({\rm INC}|{\rm DEG,PINC})$
for each combination of  DEG and PINC.
Specifically, each row corresponds to a single value of DEG, and 
each boxplot within a row corresponds  to a single value of 
PINC. The boxplot provides 
 5, 25, 50, 75 and 95\% quantiles of 
 $\Pr({\rm INC}|{\rm DEG,PINC})$. 
Note that the boxplots within a row indicate  very small increases 
in INCOME with increasing values of PINC, while 
differences across rows indicate much larger increases with DEG
(changes in the quantiles do not happen continuously  due to the binned nature of the raw data).
For high-school graduates (DEG=1), 
the estimated conditional mean incomes across levels of PINC 
are $\{ 23, 25,26,28, 29\}$ in thousands of dollars. 
For college graduates (DEG=2), the estimated means
are  $\{ 41,41,43,44,47\}$. For these mean calculations, 
 the income in a binned income category was taken as the average
 of the endpoints of the bin.

For comparison, the actual values of INC for each combination of DEG and PINC 
are plotted on the corresponding boxplots
(data are jittered to allow ties to be distinguished). 
As before, the main 
discrepancies occur for combinations of DEG and PINC 
for which there are few data. 
Also, the predictive distributions based on the copula model 
are much smoother than the empirical versions: 
The empirical conditional means of INC for DEG=1 and DEG=3 are 
 $\{ 23, 27, 24, 27,  8  \} $ and $\{ 41, 44, 35, 58, 75 \}$ 
respectively, across increasing levels of PINC. However, several of these empirical means 
are calculated from  as few as 3 or 4 samples.

\section{Notions of sufficiency}
The extended rank likelihood described above can be viewed as a
generalization of
marginal likelihood,
a standard technique
for dealing with nuisance parameters (see Section 8.3 of
\citet{severini_2000}
for a review).
One benefit of using such a likelihood is a gain in robustness,
as inference no longer depends on assumptions about the
relationship of the data to the  nuisance
parameters. Another benefit is a general simplification of the
estimation problem, as the need to estimate a potentially
high-dimensional set of parameters is eliminated.
These benefits come at the cost of potentially losing information
about the parameters of interest by only using part of the available data.
Ideally, the statistic that generates the marginal likelihood
is ``partially sufficient''
in the sense that it contains all relevant information  in the data
about
the parameter of interest.
Various definitions of partial sufficiency have been developed:
\citet{fraser_1956}
defined  $S$-sufficiency via properties of the
marginal and conditional distributions of the statistic and the data.
The concept of $G$-sufficiency was introduced in \citet{barnard_1963}
as a general principle for making inference about a parameter
of interest  when the inference problem remains invariant
under a group of transformations.
\citet{remon_1984}
developed a generalization of these notions
based on profile
likelihoods called $L$-sufficiency, which has been refined and
studied by 
\citet{barndorff-nielsen_1988,barndorff-nielsen_1999}.
The general recommendation of these authors is to base inference
for a parameter of interest on the sampling
distribution of a statistic that is sufficient in some sense.

If $F_1,\ldots, F_p$ are all continuous then there
are no ties among the data, and knowledge of $\m Z\in D$
provides a complete ordering of $\{ y_{1,j},\ldots, y_{n,j}\}$
for each $j$. This information is equivalent to the information
 contained in the ranks,
and so $\Pr (\m Z \in D|\m C)$ is equivalent to the sampling
distribution of the multivariate ranks.
Following the notation of \citet{remon_1984}
we now show that the ranks $r(\m Y)$ are a $G$-sufficient statistic
in the sense of \citet{barnard_1963}: 
Let  $\m C \in \mathcal C$ describe the copula and
$\m F = \{ F_1,\ldots,F_p\} \in \mathcal F $ the marginal
distributions, and so the parameter space is
 $\Omega= \mathcal C \times \mathcal F$ and the model space is
$\mathcal P=\{ \Pr( \cdot | \omega): \omega \in \Omega\} $,
where $\Pr(\cdot | \omega)$ is a probability measure on $\mathbb R^p$
for each $\omega\in \Omega$.
Furthermore, let   $\mathcal G$ be the group of collections of $p$
continuous strictly increasing
functions, so that $\mathcal G = \{ \m G=(G_1,\ldots ,G_p) : G_j$
 is a continuous and strictly increasing function on $\mathbb R \}$.
To each $\m G\in \mathcal G$ there corresponds a one-to-one function
on $\mathcal P$ mapping  $P(\cdot  |\omega)$ to  $P(\m G^{-1}(\cdot)|\omega)$
and the model space  is  closed under the action of $\mathcal G$.
As a result,
$\mathcal G$ induces a group  $\bar {\mathcal G}=\{ f_{\m G}: \m G \in \mathcal G\}$ on $\Omega$
defined by $P( \cdot | f_{\m G} \omega )  = P( \m G^{-1}(\cdot) | \omega)$.

If the marginals are continuous the
 orbits of $\Omega$ under $\bar {\mathcal G}$ can be put into 1-1 correspondence with
$\m C$, and $\m C$ is therefore a maximal invariant parameter.
Barnard defined a statistic $t(\m Y)$ to be $G$-sufficient if it can be
 put into
1-1 correspondence with the orbits of  $\mathbb R^p$ under $\mathcal G$.
This is the case for the ranks $r(\m Y)$ of $\m Y$, and
so $r(\m Y)$ is said to be $G$-sufficient for estimation of $\m C$.
For continuous data,
the marginal distribution of the ranks  is equal to 
the extended rank likelihood,
and so basing inference on this likelihood function can been seen
as using all available, relevant information in the $G$-sufficient sense.

A notion of sufficiency that is more directly related
to maximum likelihood estimation is $L$-sufficiency:
In the context of copula modeling, a statistic
$t(\m Y)$ is said to be $L$-sufficient for $\m C$ if
\begin{itemize}
\item[A1.]  $t(\m Y_0) = t(\m Y_1) \Rightarrow
 \sup_{\{F_1,\ldots, F_p\}\in \mc F} p(\m Y_0 | \m C,  F_1,\ldots, F_p) =
  \sup_{\{F_1,\ldots, F_p\}\in \mc F} p(\m Y_1 | \m C, F_1,\ldots, F_p) $;
\item[A2.]  $p(t(\m Y)|\m C,F_1,\ldots, F_p) = p(t(\m Y)|\m C)$.
\end{itemize}
Note that
the maximum likelihood estimate of $\m C$ and its distribution
will be a function
only of an  $L$-sufficient statistic, if one exists.
If $\mc  F$ contains only continuous marginals, then
one can show directly that the ranks $r(\m Y)$ satisfy
A1 and A2 (alternatively, \citet{remon_1984} 
shows that
a $G$-sufficient statistic is also $L$-sufficient).
Thus in the continuous case, the ranks are $G$- and $L$-sufficient,
the MLE of $\m C$ is a function of the ranks alone, and
inference for $\m C$ can be based on the
distribution of the multivariate ranks, or equivalently, the
extended rank likelihood.

If the marginals are allowed to be discontinuous then the
orbits of $\Omega$ under $\bar {\mathcal G}$
cannot be put into 1-1 correspondence
with $\m C$ and so $\m C$ is not a maximal invariant. The problem
is basically that if $F_j(\cdot)$ is a discrete CDF, then
$F_j[ G_j^{-1}(\cdot)]$ does not range over the space of all
CDF's as $\m G$ ranges over $\mathcal G$.
The ranks are no longer $L$-sufficient either:
Condition A1 holds but A2 is violated because
in the discrete case the distribution of the ranks
depends on the marginal distributions. This means
that estimation based on $\Pr(r(\m Y)|\m C ,F_1,\ldots, F_p)$
requires estimation of the nuisance parameters $F_1,\ldots, F_p$.
This may not be much of an issue if the number of levels of
each variable is low, but for moderate numbers of levels we may
wonder about the variability of the estimates due to the large number
of parameters, or the need to specify a prior distribution for
the marginals $F_1,\ldots, F_p$ in the context of Bayesian estimation.
In contrast, the extended rank likelihood based on $\Pr(\m Z \in D|\m C)$
does not depend on $F_1,\ldots, F_p$, thereby  reducing the number
of parameters to estimate  and eliminating any need for a prior
distribution  on $F_1,\ldots, F_p$. Furthermore, the extended rank likelihood
is ``sufficient'' for continuous data but can be used with mixed continuous and
discrete data.
However, the concern remains that the this likelihood
may not be making full use of the information in discrete  data about
 the copula parameters
of interest.
It would  be desirable to describe precisely any potential information
loss that results from using the rank likelihood as opposed to a full likelihood
approach. Such a description could be obtained by comparing the
curvatures of the extended rank likelihood  and full likelihood surfaces, although
the complicated parameter space and likelihood functions make description difficult except
for the simplest of cases.
A general description of the information properties of the rank likelihood in the
context of copula estimation
is a current research interest  of the author.

\section{Discussion}
This article has presented an inferential procedure for copula parameters
that can be applied to mixed continuous and discrete data.
The procedure is based on a type of marginal likelihood, called
an extended rank likelihood, which does not depend on the univariate marginal distributions
of the data. The procedure therefore allows for the estimation of
dependence parameters without the burden of having to estimate the
marginal distributions. 

The data analyzed in this paper are categorical, although
some of the variables have very large numbers of categories. 
An alternative approach to the analysis of categorical 
data is log-linear modeling. 
For categorical data, a  log-linear model
 can  potentially provide a more detailed representation of complex
dependencies and interactions than can a Gaussian copula model. 
However, if the number of categories is large and the data 
are ordinal, a copula model might be more appropriate. 
The variables AGE, INC and PCHILD in this article have 60, 21 and 19 
categories respectively. Stable log-linear analysis of these data 
would require a coarsening  of these and perhaps some of the other
variables into many fewer categories, resulting in information loss. 
In contrast, the semiparametric Gaussian copula approach taken here provides a simple dependence 
model for data having arbitrary marginal distributions, discrete or continuous.

The Gibbs sampling algorithm described in Section 3.2 is quite simple 
and performs well for the data analysis in Section 4. However, 
the fact that each  $z_{i,j}$ is being sampled 
one at a time, and from a distribution that is constrained 
by the values of $\{z_{k,j}:k\neq i\}$, 
might raise concerns that the simple Gibbs sampler might mix poorly 
in some situations. If  poor mixing occurs, 
 one remedy is to add Metropolis-Hastings 
updates that propose simultaneous changes to multiple $z_{i,j}$'s. 
One such procedure that I have implemented is to propose changes 
to the set $\{z_{i,j}: i=1,\ldots,n\}$ by shuffling the distances between 
the order statistics. 
In the examples I have tried, this type of procedure 
has given reasonable acceptance rates 
and has reduced autocorrelation. 

Inference on the scale of the  original data can be obtained 
with a posterior predictive distribution based on plugging in the 
empirical univariate marginal distributions as described in Section 4.3. 
Alternatively, a predictive distribution which accounts for uncertainty in 
the univariate marginal distributions 
can be derived as follows:
The Gibbs sampling scheme of Section 3 can be used to generate 
a joint posterior distribution for  $\m z_1,\ldots, \m z_n$ in addition 
to a new sample $\m  z_{n+1}$, for which we do not observe $\m y$-values. 
However, if  $z_{n+1,j}$   is between two other $z_j$'s having the same 
$y_j$ value, then $y_{n+1,j}$ must equal  $y_j$ as well since the 
$g_j$'s are non-decreasing. Technically, this produces a type of 
interval probability distribution for $\m y$ \citep{weichselberger_1993}, 
and for continuous data gives  univariate marginal predictive probabilities 
equivalent to the $A_n$ procedure of \citet{hill_1968}. 
For large $n$
however, this procedure is essentially equivalent to using the 
the empirical  marginal distributions. 


Although this article has focused on semiparametric estimation of
a Gaussian copula,
 the notion of  rank likelihood is equally applicable
to other copula models:
Letting $\{ p( \m u | \theta) : \theta\in \Theta \}$
denote a parametric family of copula densities and
$\{ y_{i,j} = G_j (u_{i,j}),i=1\ldots,n, \ j=1,\ldots, p \}$ be the observed data,
the extended rank likelihood for $\theta$ is given by
$\Pr(  \max \{ u_{k,j}: y_{k,j}< y_{i,j} \} < u_{i,j} < 
  \min \{ u_{k,j}: y_{i,j}< y_{k,j} \}, i=1,\ldots, n , \ j=1,\ldots, p  |\theta)$.
Given a prior distribution on $\theta$, posterior inference can be obtained
via a Markov chain Monte Carlo algorithm which iteratively  resamples
values of $\theta$ and the $u_{i,j}$'s. However, full conditional distributions
for these unknown quantities are generally hard to come by, and an MCMC sampler
based on the Metropolis-Hastings
algorithm is required  for most models.

Code to to implement the estimation strategy outlined in Section 3, written in the {\sf R} statistical computing environment,
is provided in the Appendix.
A more detailed
open-source software
package
is downloadable from
 {\sf R}-archive at the following website:
\begin{center}
\href{http://cran.r-project.org/src/contrib/Descriptions/sbgcop.html}{\tt http://cran.r-project.org/src/contrib/Descriptions/sbgcop.html}  \\
\end{center}

\appendix
\section{R-code for Gaussian copula estimation}

\begin{footnotesize}
\begin{verbatim}
# See also http://cran.r-project.org/src/contrib/Descriptions/sbgcop.html
# 
# Preconditions: Y,     an n-observations by p-variables matrix 
#                S0,    a p x p prior covariance matrix 
#                n0,    an integer hyperparameter 
#                NSCAN, an integer number of iterations 

########## helper function
rwish<-function(S0,nu){        # sample from a Wishart distribution
  sS0<-chol(S0)
  Z<-matrix(rnorm(nu*dim(S0)[1]),nu,dim(S0)[1])%*%sS0
  t(Z)%*%Z             } 

########## starting values
n<-dim(Y)[1]
p<-dim(Y)[2]
set.seed(1)
Z<-qnorm(apply(Y,2,rank,ties.method="random")/(n+1))
Zfill<-matrix(rnorm(n*p),n,p)
Z[is.na(Y)]<-Zfill[is.na(Y) ]
Z<- t( (t(Z)-apply(Z,2,mean))/apply(Z,2,sd) )
S<-cov(Z)

########## constraints 
R<-NULL
for(j in 1:p) { R<-cbind(R, match(Y[,j],sort(unique(Y[,j])))) }

########## start of Gibbs sampling scheme
for(nscan in 1:NSCAN) {

  #### update Z[,j]
  for(j in sample(1:p)) {
    Sjc<- S[j,-j]%*%solve(S[-j,-j])
    sdj<- sqrt( S[j,j] -S[j,-j]%*%solve(S[-j,-j])%*%S[-j,j]  )
    muj<- Z[,-j]%*%t(Sjc)

    for(r in sort(unique(R[,j]))){
      ir<- (1:n)[R[,j]==r & !is.na(R[,j])]
      lb<-suppressWarnings(max( Z[ R[,j]<r,j],na.rm=T))
      ub<-suppressWarnings(min( Z[ R[,j]>r,j],na.rm=T))
      Z[ir,j]<-qnorm(runif(length(ir),
               pnorm(lb,muj[ir],sdj),pnorm(ub,muj[ir],sdj)),muj[ir],sdj)
                                  }
    ir<-(1:n)[is.na(R[,j])]
    Z[ir,j]<-rnorm(length(ir),muj[ir],sdj)
                       }

  #### update S
  S<-solve(rwish(solve(S0*n0+t(Z)%*%Z),n0+n))
                      }
########## end of Gibbs sampling scheme
\end{verbatim}
\end{footnotesize}

\bibliographystyle{plainnat}
\bibliography{hoff_erlcop}

\end{document}